\documentclass{article}
\usepackage{epsfig}
\usepackage{psfrag}

\usepackage{amsmath,amssymb,amsthm}

\title{
On the continued fraction expansion of certain Engel series} 
\author{Andrew Hone\thanks{School of Mathematics,
Statistics and Actuarial Science, University of
Kent, Canterbury CT2 7NF, U.K. ~~E-mail: A.N.W.Hone@kent.ac.uk}
}

\newcommand{\beq}{\begin{equation}}
\newcommand{\eeq}{\end{equation}}
\newcommand{\bear}{\begin{array}}
\newcommand{\eear}{\end{array}}
\newcommand\la{{\lambda}}
\newcommand\La{{\Lambda}}

\newcommand\al{{\alpha}}

\newcommand\bM{\mathbf{M}}

\newtheorem{thm}{Theorem}[section]

\newtheorem{propn}[thm]{Proposition}

\newtheorem{rem}[thm]{Remark}

\newtheorem{exa}[thm]{Example}

\newtheorem{cor}[thm]{Corollary}

\newenvironment{prf}{\trivlist \item [\hskip
\labelsep {\bf Proof:}]\ignorespaces}{\qed \endtrivlist}

\theoremstyle{remark}

\newcommand{\Z}{{\mathbb Z}}

\newcommand{\seqnum}[1]{\href{http://oeis.org/#1}{\underline{#1}}}

\begin{document}

\maketitle

\begin{abstract}
An Engel series is a sum of the reciprocals of an increasing sequence of positive integers, 
which is such that each term is divisible by the previous one. 
Here we consider a particular class of Engel series, 
for which each term of the sequence is divisible by the square of the preceding one, 
and find an explicit expression for the continued fraction expansion of the sum of a generic 
series of this kind. As a special case, this includes certain series whose continued fraction expansion was found by Shallit. 
A family of examples generated by nonlinear recurrences with the 
Laurent property is considered in detail, along with some associated transcendental numbers.

\vspace{.05in}

\vspace{.05in}
\end{abstract}

\section{Introduction}

\setcounter{equation}{0}
Given a sequence of positive integers $(x_n)$, which is such that $x_n|x_{n+1}$ for all $n$, the sum of the reciprocals is 
the Engel series 
\beq \label{engel}
\sum_{j=1}^\infty \frac{1}{x_j} = \sum_{j=1}^\infty \frac{1}{y_1y_2\cdots y_j},
\eeq 
where $y_1=x_1$ and $y_{n+1}=x_{n+1}/x_n$ for $n\geq 1$. (It should be assumed that  $(x_n)$ is eventually increasing, which guarantees the convergence of the sum. A brief introduction to Engel series can be found in \cite{duve}.) In recent work \cite{curious}, we considered some particular series of this kind that are generated by certain nonlinear recurrences of second order, and are such that the sequences $(x_n)$ and $(y_n)$ appear interlaced in the continued fraction expansion. Here we start from a sequence with the stronger property that $x_n^2|x_{n+1}$. Any initial 1s can be ignored, but in what follows it will be convenient to start with $x_1=1$ and take $x_2\geq 2$, which implies $x_n\geq 2^{2^{n-2}}$ for $n\geq 2$, and we may write 
\beq\label{xdef} 
x_1=1, \qquad x_n =\prod_{j=2}^n z_j^{2^{n-j}} 
\eeq 
for some sequence of positive integers 
$(z_n)$ with $z_2\geq 2$. The corresponding Engel series is  
\beq \label{engels}
S:=\sum_{j=1}^\infty \frac{1}{x_j} =1+ \sum_{j=2}^\infty \frac{1}{z_2^{2^{j-2}} z_3^{2^{j-3}} \cdots z_j}.
\eeq 
Our main result will be to show that in the generic case, when $z_2\geq 3$ and $z_n \geq 2$ for $n\geq 3$, the continued fraction expansion of $S$ has a universal structure,  which we present explicitly.  

There is a precedent for these results in the work of Shallit, who first found the continued fraction expansion of the sum 
\beq \label{kempner}
\sum_{k=0}^\infty \frac{1}{u^{2^k}}
\eeq 
in \cite{sh1}, for integer $u\geq 3$ (with $u=2$ being a degenerate case), and went on \cite{sh3} to obtain the continued fraction for the more general series 
\beq \label{shallit}
\sum_{k=0}^\infty \frac{1}{u^{c_k}}, 
\eeq 
where $(c_k)$ is a sequence of positive integers with some non-negative $N$ such that 
$
d_n:=c_{n+1}-2c_n\geq 0 
$ 
for all $n\geq N$. If we set $z_2=u^{c_0}$ and $z_j=u^{d_{j-3}}$ 
for $j\geq 3$ in (\ref{engels}), and assume $d_n\geq 0$ for all 
$n\geq 0$, then $S-1$ coincides with (\ref{shallit}). 
 
\subsection{Outline of the paper} 

In the next section we prove the main result, namely the expression for the continued fraction expansion of a generic series of the form (\ref{engels}). Section 3 is devoted to an infinite family of examples of series of this type, which are generated by nonlinear recurrences with the Laurent property. For such nonlinear recurrence sequences, we show that the sum of the corresponding series (\ref{engels}) is a transcendental number. In the fourth section we consider the continued fractions obtained from certain degenerate cases, when either 
$z_2=2$ or $z_n =1$ for $n\geq 3$, and some particular examples of these degenerate cases are examined in more detail. The final section contains some conclusions. 

\section{Continued fractions}  
\setcounter{equation}{0}

We use the notation 
$$ 
[a_0; a_1,a_2,a_3,\ldots,a_j,\ldots ]=a_0+\cfrac{1}{a_1+\cfrac{1}{a_2+\cfrac{1}{a_3+\cdots\cfrac{1}{a_j+\cdots}}}} 
$$ 
for continued fractions, and for the 
$n$th convergent of the continued fraction $[a_0;a_1,a_2,\ldots]$ we have 
$$ 
\frac{p_n}{q_n} = [a_0;a_1,a_2,\ldots,a_n], 
$$ 
where the numerators $p_n$ and denominators $q_n$ are given in terms of the coefficients $a_j$ according to the matrix identity 
\beq\label{matid} 
\left(\begin{array}{cc} a_0 & 1 \\ 1 & 0 \end{array}\right)
\left(\begin{array}{cc} a_1 & 1 \\ 1 & 0 \end{array}\right) 
\ldots 
\left(\begin{array}{cc} a_n & 1 \\ 1 & 0 \end{array}\right) = 
\left(\begin{array}{cc} p_n & p_{n-1} \\ q_n & q_{n-1} \end{array}\right).
\eeq 
For what follows, it will also be convenient to note the identity obtained by taking the determinant of each side of (\ref{matid}), that is 
\beq\label{detid} 
p_nq_{n-1}-p_{n-1}q_n = (-1)^{n+1}.
\eeq 
For other basic results on continued fractions, the reader is referred to \cite{cassels}. 

To begin with, it is convenient to treat the factors $z_2,z_3,\ldots$ of the sequence $(x_n)$ as variables. For the first few partial sums we find the continued fraction expansions 
$$ 
S_1=1,\quad S_2=\frac{z_2+1}{z_2}=[1;z_2],\quad 
S_3=\frac{z_2^2z_3+z_2z_3+1}{z_2^2z_3}=[1;z_2-1,1,z_3-1,z_2], 
$$ 
where the 
$n$th partial sum of (\ref{engels}) is denoted $S_n$. In general it is straightforward to write $S_n$ as a fraction, that is 
\beq \label{numden}
S_n=\frac{\sum_{j=1}^{n-1}\prod_{k=2}^j z_k^{2^{n-k}-2^{j-k}}
\prod_{\ell=j+1}^n z_\ell^{2^{n-\ell}}+1} 
{z_2^{2^{n-2}}z_3^{2^{n-3}}\cdots z_n}, 
\eeq  
where the denominator is $x_n$ as given in (\ref{xdef}); 
but the continued fraction expansion of the $n$th partial sum 
is best described recursively.

The basic pattern can be seen by looking at the continued fraction for the fourth partial sum, which  is 
$$ 
S_4=[1;z_2-1,1,z_3-1,z_2,z_4-1,1,z_2-1,z_3-1,1,z_2-1].
$$ 
Observe that the first five coefficients are the same as those of $S_3$, followed by $z_4-1$, 1, and then four more coefficients which almost coincide with the last four in $S_3$ in reverse order, except that there is $z_2-1$ in place of the final $z_2$ in $S_3$.  This pattern persists, as described by the following 

\begin{propn}\label{maincf} 
Given the initial set of coefficients 
$$ 
[a_0^{(3)}; a_1^{(3)},a_2^{(3)},a_3^{(3)},a_4^{(3)}]=[1;z_2-1,1,z_3-1,z_2] 
$$ 
of the continued fraction expansion of $S_3$, of length $\ell_3=4$, define a sequence of sets of coefficients $(\{a_j^{(n)}\}_{j=0}^{\ell_{n}-1})$ with 
$\ell_n=3\cdot 2^{n-2}-1$ for $n=3,4,5,\ldots$ recursively according to 
$$ 
a_j^{(n+1)}=a_j^{(n)} \quad \mathrm{for}\quad j=0,\ldots, \ell_{n}-1,
$$ 
$$ 
a_{\ell_{n}}^{(n+1)}=z_{n+1}-1, \qquad 
a_{\ell_{n}+1}^{(n+1)}=1, \qquad a_{\ell_{n}+2}^{(n+1)}= a_{\ell_{n}-1}^{(n)}-1, 
$$ 
and 
$$ 
a_j^{(n+1)}=a_{2\ell_{n}-j+1}^{(n)} \quad \mathrm{for}\quad j=\ell_{n}+3,\ldots, 2\ell_{n}.
$$ 
Then the $n$th partial sum of the series (\ref{engels}) has the  continued fraction expansion 
\beq\label{sncf}
S_n=[a_0^{(n)}; a_1^{(n)},\ldots, a_{\ell_n -1}^{(n)}]. 
\eeq
\end{propn} 
\begin{prf} This can be done similarly 
to the proof in \cite{sh3}, but we 
prefer to use matrix computations, in the same vein as \cite{wu}. The case $n=3$ is easily verified directly. Proceeding by induction, suppose that the continued fraction expansion of $S_n$ is given by (\ref{sncf}), with  coefficients 
$a_j^{(n)}$ for $j=0,\ldots \ell_n -1$ with $\ell_n =3\cdot 2^{n-2}-1$ defined according the prescription above, and denote the numerators and denominators of the convergents by $p_j$ and $q_j$ respectively; so the final convergent gives 
$$ 
S_n = \frac{ {p}_{\ell_{n}-1}}{ {q}_{\ell_{n}-1}} , \qquad   {q}_{\ell_{n}-1}=x_n.
$$
Then for the next finite continued fraction defined by this  recursive procedure there are a total of $\ell_{n+1}=2\ell_n +1=3\cdot 2^{n-1}-1$ coefficients, that is $a_j^{(n+1)}$ for $j=0,\ldots \ell_{n+1} -1$, and for the convergents 
we use $\tilde{p}_j$,$\tilde{q}_j$ to denote numerators/denominators, respectively. So by  (\ref{matid}) we have 
\beq\label{tildes}
\bear{l}
\left(\begin{array}{cc} \tilde{p}_{\ell_{n+1}-1} &  \tilde{p}_{\ell_{n+1}-2} \\  \tilde{q}_{\ell_{n+1}-1} &  \tilde{q}_{\ell_{n+1}-2} \end{array}\right)  = 
\bM_n 
\left(\begin{array}{cc} z_{n+1}-1 & 1 \\ 1 & 0 \end{array}\right)\left(\begin{array}{cc} 1 & 1 \\ 1 & 0 \end{array}\right) 
\tilde{\bM}_n^T
,
\eear 
\eeq 
where 
$$ 
\bM_n= \left(\begin{array}{cc} a_0^{(n)} & 1 \\ 1 & 0 \end{array}\right)\cdots \left(\begin{array}{cc} a_{\ell_n-1}^{(n)} & 1 \\ 1 & 0 \end{array}\right)=
\left(\begin{array}{cc} p_{\ell_n-1} & p_{\ell_n-2} \\ q_{\ell_n-1} & q_{\ell_n-2} \end{array}\right)
$$ 
and 
$$ 
\bear{l}
\tilde{\bM}_n^T=\left(\begin{array}{cc} a_{\ell_n-1}^{(n)}-1 & 1 \\ 1 & 0 \end{array}\right)\left(\begin{array}{cc}  a_{\ell_n-2}^{(n)} & 1 \\ 1 & 0 \end{array}\right) \cdots  \left(\begin{array}{cc} a_1^{(n)} & 1 \\ 1 & 0 \end{array}\right) \\ 
    \qquad  = 
\left(\begin{array}{cc} a_{\ell_n-1}^{(n)}-1 & 1 \\ 1 & 0 \end{array}\right) \left(\begin{array}{cc} a_{\ell_n-1}^{(n)} & 1 \\ 1 & 0 \end{array}\right)^{-1} 
\bM_n^T  \left(\begin{array}{cc} a_0^{(n)} & 1 \\ 1 & 0 \end{array}\right)^{-1}\\ 
 \qquad = 
  \left(\begin{array}{cc} q_{\ell_n-1}-q_{\ell_n-2} & p_{\ell_n-1}-q_{\ell_n-1}+q_{\ell_n-2}-p_{\ell_n-2} \\ q_{\ell_n-2} & p_{\ell_n-2}-q_{\ell_n-2} \end{array}\right), 
\eear 
$$
with $T$ denoting transpose. Thus the equation  (\ref{tildes}) simplifies to yield 
$$ 
\bear{l}
\left(\begin{array}{cc} \tilde{p}_{\ell_{n+1}-1} &  \tilde{p}_{\ell_{n+1}-2} \\  \tilde{q}_{\ell_{n+1}-1} &  \tilde{q}_{\ell_{n+1}-2} \end{array}\right)  =  
\left(\begin{array}{cc} z_{n+1}q_{\ell_n-1} p_{\ell_n-1}+ 1 & 
 z_{n+1} p_{\ell_n-1}\Delta_n -1 \\   z_{n+1} q_{\ell_n-1}^2 & z_{n+1} q_{\ell_n-1}\Delta_n-1\end{array}\right), 
\eear 
$$
with $\Delta_n= p_{\ell_n-1}- q_{\ell_n-1}$, where we have used the fact that 
$$ 
\det \bM_n = \left|\begin{array}{cc} p_{\ell_n-1} & p_{\ell_n-2} \\ q_{\ell_n-1} & q_{\ell_n-2} \end{array}\right|=-1
$$ 
by  (\ref{detid}), since $\ell_n$ is odd. Hence we have 
$$ 
\tilde{p}_{\ell_{n+1}-1}= z_{n+1}q_{\ell_n-1} p_{\ell_n-1}+ 1, \qquad  \tilde{q}_{\ell_{n+1}-1}= z_{n+1} q_{\ell_n-1}^2=x_{n+1}, 
$$ 
so that 
$$ 
S_{n+1}=S_n+\frac{1}{x_{n+1}}= \frac{ {p}_{\ell_{n}-1}}{ {q}_{\ell_{n}-1}}+\frac{1}{z_{n+1}{q}_{\ell_{n}-1}^2}=\frac{\tilde{p}_{\ell_{n+1}-1}}{\tilde{q}_{\ell_{n+1}-1}}
$$ 
which is the required result. 
\end{prf} 
\begin{rem} Note that, mutatis mutandis, both the recursive structure of the partial sums and the above inductive proof hold for the partial sums of an Engel series (\ref{engel}) if, for some positive integer $n_0$, the sequence 
$(x_n)$ satisfies the weaker  condition that $x_n^2\vert x_{n+1}$ for 
$n\geq n_0$ only. This is the analogue of the fact that for the series 
(\ref{shallit}) in \cite{sh3}, 
$c_{n+1}-2c_n\geq 0$ need only hold 
for  $n\geq N$, for some $N$.
\end{rem} 
The finite continued fraction expansions of the partial sums immediately yield the continued fraction for the full series (\ref{engels}), at least for a  
generic choice of factors $z_2,z_3,\ldots$ of the sequence $(x_n)$. 
\begin{thm}\label{main}
For integer factors $z_2\geq 3$ and $z_n\geq 2$ for all $n\geq 3$, the Engel  series (\ref{engels}) has the continued fraction expansion 
\beq\label{scf} 
S=[a_0; a_1,a_2,\ldots, a_j, \ldots ] = [1; z_2-1,1,z_3-1,z_2, z_4-1,1,\ldots] 
\eeq 
where the coefficients are given by 
$
a_j = \lim_{n\to\infty}a_j^{(n)}. 
$ 
\end{thm} 
\begin{prf} The result follows from taking the limit $n\to\infty$ in (\ref{sncf}), provided that none of the coefficients in the finite continued fractions vanish. To see that this is so, note that only $1,z_2-1,z_2$ and $z_3-1$ appear as coefficients in the continued fraction for $S_3$, and all of these are non-zero with the above conditions on the $z_j$. At each step of the recursion in Proposition \ref{maincf} only  $z_{n+1}-1\geq 1$ 
and $ a_{\ell_{n}+2}^{(n+1)}= a_{\ell_{n}-1}^{(n)}-1$ are potentially new coefficients, so we must check that $ a_{\ell_{n}-1}^{(n)}-1$ cannot vanish. 
For $n=3$  the last coefficient is $a_{4}^{(3)}=z_2$, so in $S_4$ this gives  $a_{7}^{(4)}=z_2-1$, while for $n\geq4$ we have 
$ a_{\ell_{n}-1}^{(n)}=a_{1}^{(n-1)}=z_2-1$, 
so $a_{\ell_{n}+2}^{(n+1)}=z_2-2$, which is why we require  $z_2\geq 3$. Thus the only numbers that appear as 
coefficients  in the continued fraction expansion of 
$S$ are $1,z_2,z_2-2$ and $z_j-1$ for $j\geq 2$, and none of these are zero. 
\end{prf}

\section{Nonlinear recurrence sequences}  
\setcounter{equation}{0}

Among nonlinear recurrences of the form 
\beq 
\label{recf} 
x_{n+N} \, x_n = f(x_{n+1}, \ldots ,  x_{n+N-1}), 
\eeq 
where $f$ is a polynomial in $N-1$ variables, 
there is a multitude of examples which surprisingly generate integer sequences.  In a wide variety of cases,  the recurrence (\ref{recf}) has the Laurent 
property: for certain special  choices of $f$, all of the iterates belong to the ring $\Z[x_0^{\pm 1},\ldots ,x_{N-1}^{\pm 1}]$;  
as a consequence, 
if all the initial values are 1 (or $\pm 1$), then each term of the sequence is an integer. 
Such sequences were popularized by Gale \cite{gale,gale2}, and subsequently  Fomin and Zelevinsky found  a useful technique - 
 the Caterpillar Lemma \cite{fz} - which can be used to prove the   Laurent property in many cases, i.e. for recurrences coming from cluster algebras \cite{fz1} 
or in the more general setting of  Laurent Phenomenon (LP) algebras \cite{lp}.  

In \cite{honepla,numberpoly} we classified recurrences of second order, 
of the form 
\beq\label{rec2}
x_{n+2} \, x_n = f(x_{n+1}).
\eeq 
For the Laurent property to hold, 
the recurrence 
(\ref{rec2}) must belong to one of  three classes, depending on the form of $f$: (i) $f(0)\neq 0$, in which case one can apply the framework of cluster algebras (when $f$ is a binomial) or LP algebras (when it is not); 
(ii)  $f(0)= 0$, $f'(0)\neq 0$; (iii) $f(0)=f'(0)= 0$. 
In the first two classes there are additional requirements on $f$, but in the third class one can take $f(x)=x^2 F(x)$ with arbitrary $F\in\Z[x]$. 

In \cite{curious} we considered the case that $f(x)=x^2 F(x)$,  where 
$F$ has positive integer coefficients with  $F(0)=1$, and obtained the continued fraction expansion of the sum $\sum_{j=1}^\infty\frac{1}{x_j}$. In order to obtain an Engel series of the form 
(\ref{engels}), we should instead choose $F$ so that (\ref{rec2}) becomes 
\beq\label{nrec}
x_{n+2} \, x_n = x_{n+1}^{d_1} G(x_{n+1}), 
\eeq 
where 
\beq\label{reqs} 
d_1\geq 3, \quad 
G(x)\in\Z_{\geq 0}[x], \quad \deg G =d_2\geq 0, \quad G(0)\neq 0, \quad G(1)\geq 3. 
\eeq 
From  
$x_{n+2}/x_{n+1}^2 = \frac{x_{n+1}}{x_n}\cdot x_{n+1}^{d_1-3} G(x_{n+1})$,  
we see  by induction that, starting with the initial values $x_0=x_1=1$, 
$(x_n)$ is a sequence of positive integers such that $x_n^2|x_{n+1}$ 
with $x_2=G(1)\geq 3$; hence also $z_n=x_n/x_{n-1}^2\geq 3$ for $n\geq 2$. 
\begin{exa}
Taking $d_1=3$, and $G(x)=3$ for all $x$ (so $d_2=0$), the recurrence (\ref{nrec}) becomes 
$x_{n+2}x_n = 3x_{n+1}^3$, 
which generates the sequence beginning  
$
1,1,3,81,531441,5559060566555523, 
\ldots 
$.  
In this case the recurrence can be solved explicitly to yield 
$$
x_n =3^{s_n}, \qquad s_n = \frac{5-\sqrt{5}}{10}\left(\frac{3+\sqrt{5}}{2}\right)^n +  \frac{5+\sqrt{5}}{10}\left(\frac{3-\sqrt{5}}{2}\right)^n-1. 
$$     
The sum of the reciprocals is 
$$ 
S=1+\frac{1}{3}+ \frac{1}{3^4}+\frac{1}{3^{12}}+\frac{1}{3^{33}}+\cdots =[1; 2,1,8,3,80,1,2,8,1,2,19682,\ldots].
$$ 
\end{exa} 
\begin{rem} In the above example, $S-1$ is a sum of the type (\ref{shallit}) 
considered in \cite{sh3}, 
with $u=3$ and $c_n =s_{n-2}$. More generally, any recurrence of the form 
$x_{n+2}x_n = ux_{n+1}^{d_1}$ with $d_1\geq 3$ and $x_0=x_1=1$ generates a sum of this type. 
\end{rem} 
For most choices of $G$  it is not possible to give the general solution of the nonlinear recurrence (\ref{nrec}) in closed form.  
Nevertheless, one can adapt the methods 
of Aho and Sloane \cite{ahos} to write a formula 
giving precise asymptotic information. By 
rewriting (\ref{nrec}) in terms of logarithms we find that $\La_n=\log x_n$ satisfies 
\beq\label{lineq} 
\La_{n+1}-(d_1+d_2)\La_{n}+\La_{n-1} =\log c +\al_n, \quad \mathrm{with}\quad 
\al_n =\log\left(\frac{G(x_n)}{cx_n^{d_2}}\right), 
\eeq 
where 
$ G(x) = cx^{d_2}+\,\mathrm{lower}\,\, \mathrm{order}. 
$
Since $\al_n =\log(1+O(x_n^{-1}))=O(x_n^{-1})$ as $n\to\infty$, the leading order behaviour of $\La_n$ is determined by the linear expression on the left-hand side of (\ref{lineq}), 
which has the characteristic  equation $\la^2-(d_1+d_2)\la +1=0$, with 
largest root 
\beq\label{root} 
\la=\frac{d_1+d_2+\sqrt{(d_1+d_2)^2-4}}{2}>2.
\eeq 
The next two statements are equivalent to analogous formulae obtained for the sequences considered in \cite{curious}. 
\begin{propn} \label{exactf} 
For the initial conditions $x_0=x_1=1$, the logarithm $\Lambda_n=\log x_n$ of each term 
of the sequence satisfying (\ref{nrec}) is given by the formula  
\beq\label{exact} 
\La_n = \left(\frac{(1-\la^{-1})\la^n -(1-\la)\la^{-n}}{\la -\la^{-1}}-1\right) 
\log c^{-\frac{1}{d_1+d_2-2}}
+\sum_{k=1}^{n-1} \left(\frac{\la^{n-k} -\la^{k-n}}{\la -\la^{-1}}\right)\al_k, 
\eeq 
where $\al_k$ is defined as in (\ref{lineq}) and $\la$ as in (\ref{root}).  
\end{propn} 
\begin{cor} \label{cor}
To leading order, the asymptotic approximation of 
the logarithm $\La_n$ is given by
\beq\label{asymp} 
\La_n\sim C\la^n,  
\eeq 
where 
$$ 
C=  \frac{1}{d_1+d_2-2}\left(\frac{1-\la^{-1}}{\la -\la^{-1}}\right) \log c +\frac{1}{\la -\la^{-1}}\sum_{k=1}^\infty \la^{-k}\al_k, 
$$ 
and for the terms of the sequence 
$
x_n \sim c^{-\frac{1}{d_1+d_2-2}}\exp (C\la^n ).  
$  
\end{cor} 
The asymptotic behaviour of these nonlinear recurrence sequences is enough to show that the sum of the corresponding Engel series is transcendental.
\begin{thm}\label{trans} 
Suppose that the sequence $(x_n)$ with initial values $x_0=x_1=1$ is generated by the recurrence (\ref{nrec}) for some $G$ satisfying the conditions (\ref{reqs}). Then the sum $S$ in  (\ref{engels}) is a transcendental number.   
\end{thm} 
\begin{prf} This is essentially identical to the proof of Theorem 4 in \cite{curious},  so here we only sketch the argument.  Recall that 
Roth's theorem says that if $\al$ is an irrational algebraic number then for an arbitrary fixed $\delta>0$ there are only 
finitely many rational approximations $p/q$ for which 
\beq\label{roth}
\left|\al -\frac{p}{q}\right|<\frac{1}{q^{2+\delta}}. 
\eeq 
The number $S$ has an infinite continued fraction expansion, so it is irrational. 
From the asymptotics (\ref{asymp})  it follows that 
for any $\epsilon>0$ the growth condition 
\beq\label{growth}
x_{n+1}>x_n^{ \la-\epsilon}
\eeq 
 holds for all sufficiently large $n$. By making a comparison with a geometric sum, this  gives 
$$ 
\left|S-\frac{p_{\ell_n-1}}{q_{\ell_n-1}}\right|=\sum_{j=n+1}^\infty \frac{1}{x_j} <
\frac{1}{x_{n}^{\la-\epsilon-\epsilon'}}=\frac{1}{q_{\ell_n-1}^{\la-\epsilon-\epsilon'}}
$$
for any $\epsilon'>0$ and $n$  large enough. 
 So if $\epsilon$ and $\epsilon'$ are chosen such that $\la-\epsilon-\epsilon'=2+\delta >2$, then $\al=S$ has infinitely many 
rational  approximations satisfying (\ref{roth}), and hence must be transcendental. 
\end{prf}
\begin{exa}
Taking $d_1=3$, and $G(x)=2x+1$ for all $x$ (so $d_2=1$), the recurrence (\ref{nrec}) becomes 
\beq \label{nex} 
x_{n+2}x_n = x_{n+1}^3(2x_{n+1}+1), 
\eeq 
which generates the sequence 
$$ 
1,1,3,189,852910317, 
5599917937724687764238078261637795,
\ldots 
$$  
with leading order asymptotics 
$$ 
x_n \sim \frac{e^{ C (2+\sqrt{3})^n}}{\sqrt{2}}, \qquad C\approx 0.107812043.
$$ 
The sum of the reciprocals is the  transcendental number 
$$ 
S =[1; 2,1,20,3,23876,1,2,20,1,2,7697947188058154,\ldots]\approx 1.3386243.
$$ 
\end{exa} 

Engel series of the form (\ref{engels}) can also be generated by nonlinear recurrences of higher order. For instance, one can take a recurrence of third order, 
\beq\label{third} 
X_{n+3}X_n =X_{n+1}^{e_1}X_{n+2}^{e_2}H(X_{n+1},X_{n+2}),\quad\mathrm{with}\,\, e_1\geq 1,\,e_2\geq 2,
\eeq 
where the polynomial $H(X,Y)\in\Z_{\geq 0}[X,Y]$ is not divisible by either of its arguments. It is straightforward 
to show that 
the Laurent property holds for this recurrence, and from 
$X_{n+3}/X_{n+2}^2 =\frac{X_{n+1}}{X_{n}}\cdot X_{n+1}^{e_1-1}X_{n+2}^{e_2-2}H(X_{n+1},X_{n+2})$, 
it is easy to see by induction that the initial values $X_0=X_1=X_2=1$ generate an integer sequence with 
$X_n^2|X_{n+1}$ for all $n\geq 0$. Thus the sum of reciprocals starting from the index 2, that is 
\beq\label{sp}
S':=\sum_{j=2}^\infty \frac{1}{X_j},
\eeq 
 is an Engel series of the form (\ref{engels}). Note that the condition $H(1,1)\geq 3$ should be imposed, in order for Theorem \ref{main} to apply to this series. 

A particular class of recurrences of the form (\ref{third}) can be obtained by factorizing the terms of a sequence satisfying (\ref{nrec}) as 
$x_n =X_n X_{n+1}$, 
which lifts the second order recurrence to 
\beq\label{thirds} 
X_{n+3}X_n =(X_{n+1}X_{n+2})^{d_1-1}G(X_{n+1}X_{n+2}).
\eeq 
For a generic polynomial $H(X_{n+1},X_{n+2})$ on the right-hand-side of (\ref{third}), it is not immediately 
obvious which term 
will be dominant as $n\to\infty$, but in the special case (\ref{thirds}) the same techniques as for the second order recurrence can be applied directly, to show that the  leading order asymptotics is 
$
\log X_n \sim C' \la^n 
$ 
for some $C'>0$, where $\la$ is given by (\ref{root}). This means that Theorem \ref{trans} applies to the series (\ref{sp}) 
as well. 
\begin{exa}
Setting $x_n=X_nX_{n+1}$ in (\ref{nex}) gives the recurrence 
\beq \label{nexlift} 
X_{n+3}X_n = X_{n+1}^2X_{n+2}^2(2X_{n+1}X_{n+2}+1), 
\eeq 
which generates the sequence beginning  
$$ 
1,1,1,3,63,13538259, 
413636490314204194515563505,
\ldots .
$$  
To leading order, 
$
\log X_n \sim C' (2+\sqrt{3})^n$, with $C'\approx 0.0227833$.
The sum of the reciprocals in (\ref{sp}) is  the  transcendental number 
$$ 
S' =[1; 2,1,6,3,3410,1,2,6,1,2,2256800700104,\ldots] \approx  1.3492064.
$$ 
\end{exa} 
For other examples of transcendental numbers whose complete continued fraction expansion is known, see \cite{ds} and references.   

\section{Degenerate cases} 
\setcounter{equation}{0}

If either $z_2=2$ or  $z_n=1$ for some $n\geq3$, then one of the coefficients in the continued fraction 
becomes zero, and Theorem \ref{main} is no longer valid. To 
obtain a continued fraction 
with non-zero coefficients, one can use the replacement rule 
$
[\ldots, a,0,b, \ldots ] \longrightarrow [\ldots, a+b,\ldots] 
$ 
(see Proposition 3 in \cite{vdpsh}) to remove the zero. Each such replacement, decreases the length of a finite continued 
fraction  by two, so in degenerate cases the length of the continued fraction expansion of 
$S_n$ is typically shorter than the generic value $\ell_n=3\cdot 2^{n-2}-1$. Here we present the expansion for two particular 
degenerate cases, omitting details of the proof. 

\subsection{The case $z_2=2$} 

For generic values of the factors $z_j$, the sequence of lengths $\ell_n$ of partial sums begins $1,2,5,11,23$ for $n=1,2,3,4,5$. When $z_2=2$ and  $z_j\geq2$ for all $j\geq3$, the first few continued fractions for the partial sums of $S$ are $S_1=1$, 
$$ 
S_2=[1;2],\quad 
S_3=[1;1,1,z_3-1,2], \quad
S_4=[1;1,1,z_3-1,2,z_4-1,1,1,z_3-1,2], 
$$ 
which are of the same length as in the generic case, except for $S_4$ being of length 10, since at the end $[\ldots, 1,1]\to 
[\ldots, 2]$. The first zero appears in $S_5$, which contains  a single coefficient $z_2-2$, so removing this 
and making the final replacement $[\ldots, 1,1]\to [\ldots, 2]$, 
results in the length being 20: 
$$ 
S_5=[1;1,1,z_3-1,2,z_4-1,1,1,z_3-1,1,1,z_5-1,2,z_3-1,1,1,z_4-1,2,z_3-1,2].
$$ 
Thereafter the pattern continues with the continued fraction doubling in length at each step, as described by the 
following 
\begin{thm}\label{z2} 
When $z_2=2$ and $z_j\neq 1$ for $j\geq 3$, 
the  Engel series (\ref{engels}) has the continued fraction expansion 
$$
S=[a_0; a_1,a_2,\ldots, a_j, \ldots ] = [1; 1,1,z_3-1,2, z_4-1,1,\ldots] 
$$
with coefficients given by 
$
a_j = \lim_{n\to\infty}a_j^{(n)}, 
$ 
where  $(\{a_j^{(n)}\}_{j=0}^{\ell_{n}-1})$, the  sequence of sets of coefficients of 
the finite continued fractions for partial sums $S_n$, of length
$\ell_n=5\cdot 2^{n-3}$ for $n=4,5,\ldots$, is defined  by starting from 
$$ 
[a_0^{(4)}; a_1^{(4)},\ldots,a_9^{(4)}]=[1;1,1,z_3-1,2,z_4-1,1,1,z_3-1,2], 
$$ 
and obtaining subsequent coefficients according to 
$$ 
a_j^{(n+1)}=a_j^{(n)} \quad \mathrm{for}\quad j=0,\ldots, \ell_{n}-2,
$$ 
$$ 
a_{\ell_{n}-1}^{(n+1)}=1,\qquad a_{\ell_{n}}^{(n+1)}=1,\qquad
a_{\ell_{n}+1}^{(n+1)}=z_{n+1}-1, 
$$ 
$$ 
a_j^{(n+1)}=a_{2\ell_{n}-j+1}^{(n)} \quad \mathrm{for}\quad j=\ell_{n}+2,\ldots, 2\ell_{n}-2,\qquad\mathrm{and}\quad a_{2\ell_{n}-1}^{(n+1)}=2.
$$ 
\end{thm} 
In order to obtain a sequence $(x_n)$ of this degenerate type from a second order recurrence of the form (\ref{nrec}), the 
conditions (\ref{reqs}) should be modified so that $G(1)=2$, which requires that 
$G(x)=x^{d_2}+1$ for some non-negative integer $d_2$. So the recurrence becomes 
\beq\label{mrec} 
x_{n+2} \, x_n = x_{n+1}^{d_1} (x_{n+1}^{d_2}+1),
\eeq 
with $d_1\geq 3$, $d_2\geq 0$. The results of Proposition \ref{exactf}, Corollary \ref{cor} and Theorem \ref{trans} 
all apply without alteration to  sequences obtained from (\ref{mrec}). 
%
\begin{exa}
Taking $d_1=3$, and $G(x)=x+1$ for all $x$, the recurrence (\ref{nrec}) becomes 
$x_{n+2}x_n = x_{n+1}^3(x_{n+1}+1)$, 
which generates the sequence 
$$ 
1,1,2,24,172800,37150633525248000000,
\ldots 
$$  
with asymptotics 
$x_n \sim e^{C(2+\sqrt{3})^n}$, 
$C\approx 0.06224548$.  
The sum of the reciprocals is the transcendental number 
$$ 
S =[1;1,1,5, 2,299,1,1,5,1,1,1244167199,2,5,1,1,299,\ldots]\approx  1.54167245.
$$ 
\end{exa} 

\subsection{The case $z_j=1$ for $j\geq 3$} 

If we set $z_2=u$ and all other factors $z_j=1$ then $x_n=u^{2^{n-2}}$ for $n\geq 2$; the expansion (\ref{scf}) is no longer valid because each coefficient 
$z_j-1$ becomes zero for $j\geq 3$. In that case, the sum of the reciprocals is 
$S=1+\sum_{k=0}^\infty u^{-2^k}$, so that $S-1$ 
coincides with (\ref{kempner}).
The continued fractions for the partial sums were first obtained in \cite{sh1}, and a nonrecursive description was given in \cite{sh2}. 
The sequence of lengths begins 1,2,3,5,9,17, with 
$\ell_n=2^{n-2}+1$ for $n\geq 3$, and the full continued fraction is 
$$ 
S=[1;u-1,u+2,u,u,u-2,u,u+2,u,u-2,u+2,u,u-2,u,u,u+2,u,\ldots] 
$$  
for $u\neq 2$. 
The only numbers that appear as coefficients in this continued fraction 
are $1,u-2,u-1,u,u+2$.  

However, the case $u=2$ is special, since some of the coefficients in the 
above expansion become zero.  
The sequence of lengths of partial sums starts with 1,2,3,5,7,11, and 
$\ell_n=2^{n-3}+3$ for $n\geq 4$. Only the numbers 1,2,4,6 appear as coefficients in 
the continued fraction for the series, which is 
$$ 
S=[1;1,4,2,4,4,6,4,2,4,6,2,4,6,4,4,2,4,6,\ldots]. 
$$ 

The argument used to prove Theorem \ref{trans}, 
based on Roth's theorem, 
does not apply to the partial sums of the series  (\ref{kempner}), 
since the sequence $(x_n)$ does not grow fast enough: in contrast to 
(\ref{growth}), it satisfies the recurrence $x_{n+1}=x_n^2$.
However, 
a direct proof of transcendence of (\ref{kempner}), 
valid for all integers $u\geq 2$, was first given by Kempner \cite{kempner}; various alternative proofs are collected in \cite{adamc}.

\section{Conclusions} 

We have found the continued fraction expansion for an Engel series of the special type (\ref{engels}). In some cases, coming from nonlinear recurrence sequences, it has been shown that this produces transcendental numbers. We expect that the sum (\ref{engels}) should be transcendental for any choice of the factors $z_2\geq 2$ and $z_j\geq 1$ for $j\geq 3$. 
However, we do not know of a simple way to prove this  in general.

\vspace{.05in} 
\noindent 
{\bf Acknowledgments:} 
This work is supported by Fellowship EP/M004333/1 
from the EPSRC.
The original inspiration came from Paul Hanna's  observations concerning the nonlinear recurrence sequences 
described in  \cite{curious}, which were communicated via the Seqfan mailing list. 
The author is grateful to Jeffrey Shallit for pointing out some of his work, as well as Kempner's proof \cite{kempner}.


\small
\begin{thebibliography}{99}

\bibitem{adamc}
B.~Adamczewski, 
The many faces of the Kempner number, 
{\it J. Integer Sequences} {\bf 16} (2013), Article 13.2.15.
\bibitem{ahos} A.~V.~Aho and N.~J.~A.~Sloane,
Some doubly exponential sequences, 
{\it Fibonacci Quart.} \textbf{11} (1973), 429--437.


\bibitem{cassels} J.~W.~S.~Cassels, {\it An Introduction to Diophantine Approximation}, 
Cambridge University Press, 1957. 


\bibitem{ds} J.~L.~Davison and J.~O.~Shallit, 
Continued fractions for some alternating series, 
{\it Monatsh. Math.} \textbf{111} (1991), 119--126.

\bibitem{duve} D.~Duverney, 
{\it Number Theory: An Elementary Introduction Through Diophantine Problems}, 
World Scientific, 2010. 

\bibitem{fz1} S.~Fomin and A.~Zelevinsky,
Cluster algebras I: Foundations,
{\it J. Amer. Math. Soc.} \textbf{15} (2002), 497--529.

\bibitem{fz} S.~Fomin and A.~Zelevinsky,
The Laurent phenomenon,
{\it Adv. Appl. Math.} \textbf{28} (2002),  119--144.





\bibitem{gale} D.~Gale, 
The strange and surprising saga of the Somos sequences, 
{\it Math. Intelligencer} \textbf{13 (1)} (1991),
40--42.

\bibitem{gale2}  D.~Gale, 
Somos sequence update,
{\it Math. Intelligencer} \textbf{13 (4)} (1991),
49--50.






\bibitem{honepla} A.~N.~W.~Hone,
Singularity confinement for maps with the Laurent property,
{\it Phys. Lett. A} {\bf 361} (2007),
341--345.
 


\bibitem{numberpoly} A.~N.~W.~Hone,
Nonlinear recurrence sequences and Laurent polynomials,  in 
J. McKee and C. Smyth, eds., 
\textit{Number Theory and Polynomials},
 LMS Lecture Notes Series, vol.
{\bf 352},  Cambridge, 2008, pp. 188--210.

\bibitem{curious} A.~N.~W.~Hone,
Curious continued fractions, nonlinear recurrences and transcendental numbers, 
{\it J. Integer Sequences} {\bf 18} (2015), Article 15.8.4.

\bibitem{kempner}
A.~J.~Kempner,
On Transcendental Numbers,
{\it Trans. Amer. Math. Soc.}
{\bf 17} (1916), 476--482. 
 
\bibitem{lp} T.~Lam and P.~Pylyavskyy, 
Laurent phenomenon algebras, preprint, 
{\tt http://arxiv.org/abs/1206.2611v2}. 

\bibitem{vdpsh} 
A.~J.~van der Poorten and J.~O.~Shallit,  
Folded continued fractions, 
{\it J.Number Theory} {\bf 40} (1992), 
237--250.

\bibitem{sh1} J.~O.~Shallit,  
Simple continued fractions for some irrational numbers,
{\it J. Number Theory} {\bf 11} (1979), 
209--217.

\bibitem{sh2} J.~O.~Shallit,  
Explicit descriptions of some continued fractions,
{\it Fibonacci Quart.} {\bf 20} (1982), 
77--81.


\bibitem{sh3} J.~O.~Shallit,  
Simple continued fractions for some irrational numbers. II, 
{\it  J. Number Theory} {\bf  14} (1982), 
228--231.

\bibitem{wu} T.~Wu,  
On the proof of continued fraction expansions for irrationals, 
{\it  J. Number Theory} {\bf  23} (1986),  
55--59.



\end{thebibliography}
\end{document}